\documentclass[12pt,reqno]{amsart}

\usepackage[arrow,matrix,curve]{xy}

\usepackage[dvips]{graphicx} 

\usepackage{amssymb, latexsym, amsmath, amscd, array, hyperref
}

\theoremstyle{definition}

\numberwithin{equation}{section}

\newcommand\N {{\mathbb N}} 

\newcommand\R {{\mathbb R}}

\newcommand\Los{{\L}o{\'s}}

\author[P. B\l aszczyk]
{Piotr B\l{}aszczyk}\address{P. B\l{}aszczyk, Institute of
Mathematics, Pedagogical University of Cracow,
Poland}\email{pb@up.krakow.pl}

\author
[V. Kanovei]
{Vladimir Kanovei} \address{V. Kanovei, IPPI, Moscow, and MIIT,
Moscow, Russia}\email{kanovei@googlemail.com}

\author 
[M. Katz] 
{Mikhail G. Katz}\address{M. Katz, Department of Mathematics, Bar Ilan
University, Ramat Gan 52900 Israel}\email{katzmik@macs.biu.ac.il}

\author
[D. Sherry]
{David Sherry}\address{D. Sherry, Department of Philosophy, Northern
Arizona University, Flagstaff, AZ 86011, US}
\email{David.Sherry@nau.edu}

\begin{document}

\thispagestyle{empty}


\title[Controversies in the foundations of analysis] {Controversies in
the foundations of analysis: Comments on Schubring's \emph{Conflicts}}

\begin{abstract}
\emph{Foundations of Science} recently published a rebuttal to a
portion of our essay it published two years ago.  The author,
G. Schubring, argues that our 2013 text treated unfairly his 2005
book, \emph{Conflicts between generalization, rigor, and intuition}.
He further argues that our attempt to show that Cauchy is part of a
long infinitesimalist tradition confuses \emph{text} with
\emph{context} and thereby misunderstands the significance of Cauchy's
use of infinitesimals.  Here we defend our original analysis of
various misconceptions and misinterpretations concerning the history
of infinitesimals and, in particular, the role of infinitesimals in
Cauchy's mathematics.  We show that Schubring misinterprets Proclus,
Leibniz, and Klein on non-Archimedean issues, ignores the Jesuit
\emph{context} of Moigno's flawed critique of infinitesimals, and
misrepresents, to the point of caricature, the pioneering Cauchy
scholarship of D.~Laugwitz.

\end{abstract}

\maketitle

\tableofcontents

\section{Introduction}

In our article \cite{BKS} in the \emph{Foundations of Science}, we
sought to place Cauchy's work within an uninterrupted line of an
infinitesimalist tradition, and to vindicate his infinitesimal
definitions of concepts like continuity, convergence, and the Dirac
delta function (see also \cite{KT}).  We also sought to champion the
pioneering work of D.~Laugwitz from the 1980s that challenges the
received views on Cauchy as a precursor of the Weierstrassian
\emph{epsilontik}.  We developed an approach to the history of
analysis as evolving along separate, and sometimes competing, tracks.
These are the A-track, based upon an Archimedean continuum; and
B-track, based upon what we called a Bernoullian (i.e.,
infinitesimal-enriched) continuum.

\subsection{Punctiform and non-punctiform continua}

Historians often view the work in analysis from the 17th to the middle
of the 19th century as rooted in a background notion of continuum that
is not punctiform.  This necessarily creates a tension with modern,
punctiform theories of the continuum, be it the A-type set-theoretic
continuum as developed by Cantor, Dedekind, Weierstrass, and others,
or a B-type continuum as developed by Hewitt, \Los, Robinson, and
others.  How can one escape a trap of presentism in interpreting the
past from the viewpoint of set-theoretic foundations commonly accepted
today, whether of type A or B?  

In analyzing Cauchy's work one must be careful to distinguish between
its \emph{syntactic} aspects, i.e., procedures and inferential moves,
on the one hand, and \emph{semantic} aspects related to the actual
construction of the entities such as points of the continuum, i.e.,
issues of the \emph{ontology} of mathematical entities such as points.
We found that in his work on Cauchy, Laugwitz was careful not to
attribute modern set-theoretic constructions to work dating from
before the heroic 1870s, and focused instead of Cauchy's procedures
(see also Section~\ref{s55}).

\subsection{Rebuttals}

In keeping with its tradition of fostering informed debate,
\emph{Foundations of Science} recently published a response to a
portion of our essay.  The author, G.~Schubring, argues that our
attempt to show that Cauchy is part of a long infinitesimalist
tradition confuses \emph{text} with \emph{context} and thereby
misunderstands the significance of Cauchy's use of infinitesimals.
Schubring urges, instead, that Cauchy's intentions must be
``reconstructed according to the contemporary conceptual horizon of
the related conceptual fields'' \cite[Section~3]{Sc15}.

He further insinuates that our analysis is ``contrary to a style
appropriate for a reasoned scientific discussion'' (ibid., Section~5).

We grant that we dealt too briefly with Schubring's treatment of
Cauchy.  Here we develop more fully our thesis that the book
\emph{Conflicts} \cite{Sc05}
\begin{enumerate}
\item
misinterprets both Leibniz and Felix Klein on infinitesimals;
\item
escapes the need of a genuine hermeneutical effort by ignoring the
historical \emph{context} of Moigno's flawed critique of
infinitesimals;
\item
underestimates the role of infinitesimals in Cauchy's work;
\item
 misrepresents Laugwitz's Cauchy scholarship, e.g., by alleging that
 Laugwitz tried to ``prove'' that Cauchy used hyper-real numbers; and
\item conflates Laugwitz's analysis with D.~Spalt's.
\end{enumerate}

The issue of modern interpretation of historical mathematics is a deep
and important one.  Much recent scholarship has argued that modern
conceptions shed little light on historical mathematics.  For example,
scholars have criticized late 19th and early 20th century attempts to
attribute a geometrical algebra to Euclid's Elements (e.g., Szabo
1969, Unguru 1975 and 1979).  More recently the nonsequiturs and
unwarranted assumptions of these critiques have been brought to light
\cite{Bl16}, and the present article attempts to accomplish something
similar.  It continues a program of re-evaluation of the history of
mathematical analysis undertaken in \cite{BKS}, \cite{B11},
\cite{Ba14}, and other texts.  Thus, we highlight the significance of
Simon Stevin's approach to undending decimals in \cite{KK12b}.  The
text \cite{KSS13} re-evaluates Fermat's contribution to the genesis of
the infinitesimal calculus via his technique of \emph{adequality}.
The articles \cite{KS2}, \cite{KS1}, and \cite{SK} argue that
Leibniz's theoretical strategy for dealing with infinitesimals was
more soundly based than critiques by George Berkeley and others.
Euler's infinitesimal analysis is defended against A-track re-writings
in the articles \cite{KKKS} and \cite{Ba16} .  Cauchy's infinitesimal
legacy is championed in the articles \cite{KK11}, \cite{BK},
\cite{KK12a}, \cite{TK}.

\section{Conceptual horizons}

Schubring has repeatedly emphasized the importance of context and
external history (as opposed to what he calls \emph{immanent
history}).  We may agree with Schubring that
\begin{quote}
The meaning [of Cauchy's infinitesimals] has to be reconstructed
according to the contemporary conceptual horizon of the related
conceptual fields. Thus, a genuine hermeneutical effort is
inescapable.  \cite[Section~3]{Sc15}
\end{quote}
Indeed, Cauchy's \emph{contemporary conceptual horizon} needs to be
taken into account.  However, it is equally true that the meaning of
Cauchy's infinitesimals cannot be reconstructed, as Schubring attempts
to do, while at the same time ignoring the context of the broader
\emph{conceptual horizon} that surely incorporated an awareness of
Archimedean and non-Archimedean issues stretching all the way back to
Euclid.  

We will show that Schubring misinterprets both Leibniz and Felix Klein
on both the Archimedean property and infinitesimals, and ends up
endorsing Moigno's flawed critique of infinitesimals; see
Sections~\ref{s4} and \ref{s5}.  Moreover, Schubring himself ignores
the contemporary historical \emph{context} of Moigno's critique; see
Section~\ref{ban}.

\section{Axioms for Archimedean quantities}
\label{two}

The notion of infinitesimal is related to the axiom of Archimedes,
i.e., Euclid's \emph{Elements}, Definition V.4 (see
\cite[Section~II.3]{De15} for more details).  The theory of magnitudes
as developed in Book V of the \emph{Elements} can be formalized as an
ordered additive semigroup $M$ with a total order, characterized by
the five axioms given below.

\cite{Bec} and \cite[pp.~101-122]{BM} provide detailed sources for
the axioms below in the primary source (Euclid).  See also
\cite[pp.~118-148]{Mue} which mostly follows Beckmann's development.
Axiom E1 below interprets Euclid~V.4:

\begin{enumerate}
\item[E1] $(\forall x, y\in M)(\exists n\in \N)(nx > y)$,
\item[E2] $(\forall x, y\in M)(\exists z\in M)(x < y \Rightarrow x + z=y)$,
\item[E3] $(\forall x, y, z\in M)(x < y \Rightarrow x + z < y + z)$,
\item[E4] $(\forall x\in M)(\forall n \in \N)(\exists y\in M)(x =ny)$,
\item[E5] $(\forall x, y, z\in M)(\exists v\in M)(x : y :: z : v).$
\end{enumerate}

\section{Schubring's pseudohistory}

We find multiple problems with Schubring's analysis of authors ranging
from Proclus and Leibniz to Moigno and Klein.
%
%
Here we present several examples.

\subsection{Klein's horn-angles misunderstood by Schubring}
\label{s4}

It is well known that infinitesimals form a \emph{proper} subset of a
non-Archimedean structure.  Now to formulate what it might mean for a
structure to be \emph{NON-Archimedean}, one needs first to make up
one's mind as to what \emph{Archimedean} itself might mean.
Surprisingly, the book \cite{Sc05} provides \emph{no formulation
whatsoever} concerning the so-called Archimedean axiom%
\footnote{On page 17, Schubring cites Proclus' correct reading of
Euclid V.4, but as soon as Schubring attempts to paraphrase this in
his own terms, he immediately gets it wrong by describing
infinitesimals in terms of \emph{incommensurability}.  The term
\emph{incommensurability} is used to describe phenomena related to
irrationality, including both previous occurrences of the term in
\cite{Sc05} on pages 12 and 16.}
(the term was introduced in \cite{St83}).  More importantly,
Schubring's remarks on Klein's horn-angles reveal a misconception on
Schubring's part concerning non-Archimedean systems:
\begin{quote}
[Felix] Klein has shown [sic] in detail that the hornlike angles form
a model of non-Archimedean quantities.  \cite[p.~17]{Sc05}
\end{quote}
But is that really what Klein had shown?  In Klein's approach,
horn-angles form an ordered structure isomorphic to~$\R_+$ (when
parametrized by the curvature), which can hardly be said to be
non-Archimedean.  Klein goes on to form a broader structure
incorporating both horn-angles and rectilinear angles.  That structure
is non-Archimedean (see Section~\ref{two} for a discussion of the
term), and a horn-angle constitutes an infinitesimal within this
structure.  More precisely, a horn-angle is an infinitesimal
\emph{relative} to rectilinear angles.  Thus, it is incorrect to
assert, as Schubring does, that ``hornlike angles form a model of
non-Archimedean quantities.''  Klein writes:
\begin{quote}
Thus multiplication of the angle of a tangent circle [i.e.,
horn-angle] by an integer [$n$] always yields another angle of a
tangent circle, and every multiple $n$ is necessarily smaller, by our
definition, than, say, the angle of a fixed intersecting line, however
large we take~$n$. Thus the axiom of Archimedes is not satisfied; and
the angles of tangent circles must be looked upon, accordingly, as
actually infinitely small \emph{with respect to} the angle of an
intersecting straight line.  \cite[p.~205]{Kl25} (emphasis added)
\end{quote}
Thus, Klein's horn-angles are only infinitely small \emph{with respect
to} rectilinear angles.  Schubring's imperfect understanding of even
the basic issues involving infinitesimals leads to errors of judgment
on his part with regard to both Leibniz and Cauchy, as we show in
Sections~\ref{s5} and~\ref{43}.

\subsection
{Leibniz's infinitesimals misinterpreted by Schubring}
\label{s5}

A quotation from Leibniz on infinitesimals in \cite{Sc05} is
accompanied by the following claim:
\begin{quote}
[Leibniz] declared two homogeneous quantities to be equal that
differed only by a quantity that was arbitrarily smaller than a
\emph{finite quantity}.  \cite[p.~169]{Sc05} (emphasis added)
\end{quote}
Schubring takes two homogeneous quantities $x,y$ to be equal if their
difference $x-y$ is a quantity that is arbitrarily smaller than a
finite quantity, or in symbols
\[
x=y \iff x-y< \text{each finite quantity}.
\]
However, an infinitesimal itself is a finite quantity (i.e., it is not
infinite).  How can it be smaller than \emph{each finite quantity}?
In other words, how can a quantity be smaller than its own half, its
own quarter, etc.?  There seems to be a logical inconsistency involved
in Schubring's summary of the definition of an infinitesimal.  The
inconsistency is identical to the one claimed by Moigno and endorsed
by Schubring; see Section~\ref{43}.

In fact, Leibniz is guiltless here, and it is Schubring who
misinterprets Leibniz's discussion of infinitesimals.  Leibniz wrote:
\begin{quote}
I agree with Euclid Book V Definition 5 that only those homogeneous
quantities are comparable, of which the one can become larger than the
other if multiplied by a number, that is, a \emph{finite number}.  I
assert that entities, whose difference is not such a quantity, are
equal.  [\ldots] This is precisely what is meant by saying that the
difference is smaller than any given quantity.%
\footnote{Leibniz uses the term \emph{finite} in his paraphrase of
Euclid's definition V.5 (or V.4, as discussed in footnote~\ref{V4}),
but here he is dealing with a finite \emph{integer} $n$ (which, in
modern terminology, is tending to infinity), so that $n\epsilon$
always stays less than $1$ thereby violating the Archimedean property,
if $\epsilon$ is infinitesimal.}
(Leibniz as quoted in \cite[p.~169]{Sc05})
\end{quote}
This passages deals both with comparable quantities and
infinitesimals, but \emph{finite numbers} appear only in Leibniz's
characterisation of comparable quantities, not of infinitesimals.  We
interpret Leibniz's reference to Euclid V.5 as a reference to the
Archimedean axiom%
\footnote{\label{V4}Leibniz lists number V.5 for Euclid's definition
instead of V.4.  In some editions of the \emph{Elements} this
definition does appear as V.5.  Thus, \cite{Eu} as translated by
Barrow in 1660 provides the following definition in V.V (the notation
``V.V'' is from Barrow's translation): \emph{Those numbers are said to
have a ratio betwixt them, which being multiplied may exceed one the
other}.  For our interpretation of this, see Section~\ref{two},
Axiom~E1.}
(see Section~\ref{two}).

The way Schubring uses \emph{finite quantity} in his erroneous
paraphrase of Leibniz has the effect of endowing infinitesimals with
contradictory properties.  Schubring's erroneous paraphrase is not an
isolated incident, for it is consistent with his endorsement of
Moigno's flawed critique of infinitesimals analyzed in
Section~\ref{43}.

\subsection{Schubring endorsing Moigno}
\label{43}

It is no accident that Schubring gives a favorable evaluation of
Moigno's critique of the infinitely small.  Moigno wrote in the
introduction to his 1840 book that infinitesimals are contradictory,
arguing that it is impossible for something be less than its own half,
its own quarter, etc.~(see Section~\ref{ban}) Schubring endorses
Moigno's critique of infinitesimals in the following terms:
\begin{quote}
[Moigno] not only fails to use [the \emph{infiniment petits}] as a
basic concept, but also even explains explicitly why they are
inappropriate as such.  This makes Moigno the first writer to pick
apart the traditional claim in favor of their purported
\emph{simplicit\'e} \cite[p.~445]{Sc05}
\end{quote}
Attributing such alleged \emph{picking apart} to Moigno involves a
fundamental misunderstanding on Schubring's part, closely related to
his misreading of Leibniz (see Section~\ref{s5}).  Without properly
understanding infinitesimals first, a scholar can't properly evaluate
historical criticisms of infinitesimals, either.  Schubring similarly
ignores an important aspect of the relevant \emph{external} history,
namely the Jesuit context of Moigno's remarks on infinitesimals (see
Section~\ref{ban}).

Given Schubring's endorsement of Moigno's claim that infinitesimals
are self-contradictory, it is no wonder that he should seek to save
Cauchy's reputation by attempting to minimize the significance of
infinitesimals in Cauchy's work.

\section{Schubring vs Laugwitz}

Schubring's analysis of Laugwitz's Cauchy scholarship contains
numerous inaccuracies and misrepresentations.

\subsection{Were Cauchy's theorems always correct?}
\label{s51}

Schubring comments as follows on Laugwitz's work on Cauchy's sum
theorem (a series of continuous functions under suitable conditions
converges to a continuous function):
\begin{quote}
[Giusti's 1984 article] spurred Laugwitz to even more detailed
attempts to banish the error and confirm that Cauchy had used
hyper-real numbers.  On this basis, he claims, the errors vanish and
the theorems become correct, or, rather, \emph{they always were
correct} (see Laugwitz 1990, 21).  \cite[p.~432]{Sc05} (emphasis
added)
\end{quote}
Schubring is making two separate claims here:
\begin{enumerate}
\item
that Laugwitz sought to show that Cauchy used hyper-real numbers; and
\item
that Laugwitz asserted that Cauchy's sum theorem was always correct
(including its 1821 version in the \emph{Cours d'Analyse}).  
\end{enumerate}
Claim (1) will be examined in Section~\ref{s10}.  Here we will examine
claim~(2).  Did Laugwitz assert that the theorems were always correct?
Let us consider the relevant passage from Laugwitz's text in
\emph{Historia Mathematica}:
\begin{quote}
Rather late in his life Cauchy [1853] admitted that the statement of
his theorem (but not its proof) was \emph{incorrect}: ``Au reste, il
est facile de voir comment on doit modifier l'\'enonc\'e du
th\'eor\`eme, pour qu'il n'y [ait] plus lieu \`a aucune exception''
[Cauchy 1853, 31-32].  \cite[p.~265]{La87} (emphasis added)
\end{quote}
Thus, Laugwitz acknowledges that the 1821 formulation of the sum
theorem was \emph{incorrect} as stated, and moreover that Cauchy
himself had recognized its incorrectness.  This is contrary to the
claim concerning Laugwitz's position made by Schubring.  Schubring's
claim amounts to a misrepresentation - indeed a caricature - of
Laugwitz's position.
Misrepresenting another scholar's work is, in Schubring's phrase,
``contrary to a style appropriate for a reasoned scientific
discussion.'' \cite[Section~5 `The Climax']{Sc15} In addition,
Schubring frequently conflates Laugwitz's position with that of
D.~Spalt.%
\footnote{Schubring repeats the performance in 2015 when he claims:
``I am analysing at length the methodological approach of Laugwitz
(and Spalt), which consists in attributing to Cauchy (his) own
`universe of discourse.'\,'' \cite[Section~3]{Sc15}.  But Spalt's
approach is not identical to Laugwitz's!}
Spalt's analysis was already criticized in \cite{KK11}.

\subsection{Didactic component}
\label{s6}

Schubring makes the following three successive claims concerning
Laugwitz's Cauchy scholarship:
\begin{enumerate}
\item
the controversy basically centers on such fundamental concepts as
continuity and convergence;
\item
Laugwitz himself talks about the `essentially didactic components of
the infinitesimals in Cauchy' (ibid., 18) (cf. \cite [p. 18] {La90});
\item
it is in textbooks that we find these basic terms and not in
isolated research memoirs
\end{enumerate}
(see \cite[p.~433]{Sc05}).  Schubring's claim (3) concerning `basic
terms' apparently refers to the terms \emph{continuity},
\emph{convergence}, and \emph{infinitesimal} that he mentioned earlier
in (1) and (2).  In his claim~(3), Schubring appears to assert that
these terms are found only in Cauchy's textbooks rather than in
Cauchy's research memoirs.  Schubring claims that Laugwitz describes
Cauchy's infinitesimals as `essentially didactic.'  Schubring does not
cite the relevant sentence from Laugwitz, who states:
\begin{quote}
Das Beispiel gibt Gelegenheit, auf die wesentlich didaktische
Komponente des Infinitesimalen bei Cauchy hinzuweisen.
\cite[p.~18]{La90}
\end{quote}
This can be translated as follows:
\begin{quote}
The example affords us the opportunity to point out the essentially
didactic component of Cauchy's infinitesimals.
\end{quote}
Laugwitz's use of the singular \emph{Komponente} suggests that, among
\emph{other} components (i.e., aspects) of infinitesimals, there is
\emph{also} a didactic component.  Schubring makes it appear as if
Laugwitz views Cauchy's infinitesimals as being \emph{limited} to
their didactic role (this interpretation is motivated by Schubring's
contention that Cauchy was forced by a curricular committee to include
infinitesimals in his textbook).  Schubring further reveals that he
himself adheres to such a view when he claims in (3) that Cauchy's
infinitesimals are found in textbooks but not in research memoirs.  We
will now show that both of Schubring's positions are untenable.

Neither Schubring's view of Cauchy's infinitesimals as being limited
to their didactic role, nor Schub\-ring's interpretation of Laugwitz's
comment stand up to scrutiny.  Note that Schubring's claim (3) amounts
to the extraordinary contention that Cauchy only used infinitesimals
in textbooks only but not in research articles, something Laugwitz
never claimed.

Was Cauchy's use of infinitesimals indeed limited to his textbooks, as
Schubring claims?  Certainly not.  The article \cite{Ca53} on the
corrected version of the sum theorem is a research article.  It deals
with a property closely related to uniform convergence of series of
functions, which was the cutting edge of research in analysis at the
time.

The 1853 article happens to deal with all three notions that Schubring
mentions in this paragraph, namely, \emph{continuity, convergence},
and \emph{infinitesimal}.  Indeed, continuity and convergence are
mentioned already in Cauchy's title, while infinitesimals are
exploited in the definition of continuity given already at the
beginning of the article:
\begin{quote}
\ldots{} une fonction~$u$ de la variable r\'eelle~$x$ sera
\emph{continue}, entre deux limites donn\'ees de~$x$, si, cette
fonction admettant pour chaque valeur interm\'ediaire de~$x$ une
valeur unique et finie, un accroissement \emph{infiniment petit}
attribu\'e \`a la variable produit toujours, entre les limites dont il
s'agit, un accroissement \emph{infiniment petit} de la fonction
elle-m\^eme.%
\footnote{Translation: ``A function $u$ of a real variable $x$ will be
\emph{continuous} between two given bounds on $x$ if this function,
taking for each intermediate value of $x$ a unique finite value, an
infinitely small increment given to the variable always produces,
between the bounds in question, an infinitely small increment of the
function itself.''}
\cite{Ca53} [emphasis on \emph{infiniment petit} added]
\end{quote}
Furthermore, the 1853 version of the sum theorem cannot even be
formulated without infinitesimals, and Cauchy's proof procedure uses
them in an essential way, as analyzed in \cite[pp.~264-266]{La87}.

Cauchy wrote many other research articles exploiting infinitesimals in
an essential way.  An example is his memoir \cite{Ca32}, where he
expresses the length of a curve in terms of the average of its
projections to a variable axis.  His proof procedure involves
decomposing the curve into infinitesimal portions and proving the
result for each portion.  This article is considered by specialists in
geometric probability to be a foundational text in that field; see
e.g., \cite{Hy12}.%
\footnote{Another example of Cauchy's use of infinitesimals in
research is his foundational text on elasticity \cite{Ca23} where ``un
\'el\'ement infiniment petit'' is exploited on page 302.  The article
is mentioned in \cite[p.~378]{Fr71}.}

Thus Schubring's claim (3) to the effect that Cauchy's infinitesimals
appear only in textbooks is not merely inaccurate, revealing an
incomplete knowledge of the original documents on his part, but more
importantly it obscures the fundamental role of infinitesimals in
Cauchy's thinking.

After mentioning the \emph{essentially didactic component of Cauchy's
infinitesimals} in \cite[p.~18]{La90}, Laugwitz goes on to discuss the
corrected sum theorem three pages later, starting at the bottom of
page~21.  Laugwitz was well aware that Cauchy's infinitesimals are not
limited to their didactic function, contrary to Schubring's claim~(2).
Three years earlier, Laugwitz had discussed Cauchy's corrected sum
theorem in detail in his article in \emph{Historia Mathematica}
\cite[pp.~264-266]{La87}.

In conclusion, Schubring both misrepresents Laugwitz's position on
Cauchy's infinitesimals, and reveals his (Schubring's) own
misconceptions regarding the latter.

\section{Cauchy, Laugwitz, and hyperreal numbers}
\label{s10}

In this section we examine the relationship between Cauchy's
infinitesimals and modern infinitesimals as seen by Laugwitz, and the
related comments by Fraser, Grabiner, and Schubring.

\subsection{Schubring--Grabiner spin on Cauchy and hyperreals}

Schubring's opposition to Laugwitz's interpretation of Cauchy found
expression in the following comment, already quoted in
Section~\ref{s51}:
\begin{quote}
[Giusti's article] spurred Laugwitz to even more detailed attempts to
banish the error and \emph{confirm that Cauchy had used hyper-real
numbers}.  On this basis, he claims, the errors vanish and the
theorems become correct, or, rather, they always were correct (see
Laugwitz 1990, 21).  \cite[p.~432]{Sc05} (emphasis added)
\end{quote}
The matter of Schubring's misrepresentation of Laugwitz's position
with the regard to the correctness of the statement of Cauchy's
theorems was already dealt with in Section~\ref{s51}.  In this
section, we will examine Schubring's contention that Laugwitz claimed
that Cauchy used the hyperreals.  In this passage, Schubring is
\emph{referring to} the article \cite{La90}, but he is most decidedly
not \emph{quoting} it.  In fact, there is no mention of the hyperreals
on page 21 in \cite{La90}, contrary to Schubring's claim.  What we do
find there is the following comment:
\begin{quote}
The ``mistakes'' show rather, as \emph{experimenta crucis} that one
must understand Cauchy's terms/definitions [\emph{Begriffe}], in the
spirit of the motto,%
\footnote{Here Laugwitz is referring to Cauchy's motto to the effect
that ``Mon but principal a \'et\'e de concilier la rigueur, dont je
m'\'etais fait une loi dans mon \emph{Cours d'analyse} avec la
simplicit\'e que produit la consideration directe des quantit\'es
infiniment petites.''}
in an infinitesimal-mathematical sense.  \cite[p.~21]{La90}
(translation ours)
\end{quote}
We fully endorse Laugwitz's comment to the effect that Cauchy's
procedures must be understood in the sense of infinitesimal
mathematics, rather than paraphrased to fit the \emph{epsilontik}
mode.  Note that we are dealing with an author, namely Laugwitz, who
published Cauchy studies in the leading periodicals \emph{Historia
Mathematica} \cite{La87} and \emph{Archive for History of Exact
Sciences} \cite{La89}.%
\footnote{The fact that Laugwitz had published articles in leading
periodicals does not mean that he could not have said something wrong.
However, it does suggest the existence of a strawman aspect of
Schubring's claims against him.}
The idea that Laugwitz would countenance a claim that Cauchy ``had
used hyper-real numbers'' whereas both the term \emph{hyper-real} and
the relevant construction were not introduced by E. Hewitt until 1948
\cite[p.~74]{He48}, strikes us as far-fetched.  Meanwhile, in a
colorful us-against-``them'' circle-the-wagons passage, J.~Grabiner
opines that
\begin{quote}
[Schubring] effectively rebuts the partisans of nonstandard analysis
who wish to make Cauchy one of \emph{them}, using the work of Cauchy's
disciple the Abb\'e Moigno to argue for Cauchy's own intentions.
\cite[p.~415]{Gr06}. (emphasis added)
\end{quote}
Grabiner's comment betrays insufficient attention to the
procedure/ontology distinction, while her endorsement of Moigno's
critique of infinitesimals is as preposterous as Schubring's.
Schubring indeed uses the work of Moigno in an apparent attempt to
refute infinitesimals (see Section~\ref{43}).  Moigno's confusion on
the issue of infinitesimals is dealt with in Section~\ref{ban}.

\newcommand\bash{ 
The Oxford English Dictionary defines the verb
\emph{to demonize} as follows:
\begin{quote}
To portray (a person or thing) as wicked and threatening, (now)
esp.~in an inaccurate or misrepresentative way. (Now the usual
sense.)
\end{quote}
Schubring's comments on Laugwitz misrepresent Laugwitz's position and
portray him as bent on having his way with regard to Cauchy even at
the expense of attributing to Cauchy things (``hyper-reals'') Cauchy
could not possibly have done.  This is as close as a historian can get
to portraying another historian as \emph{wicked}, in line with the OED
definition.  Our conclusion is that Schubring did not properly analyze
Laugwitz but rather misrepresented him via a strawman reading and
sought to ostracize him through ridicule and demonisation.}
%

Contrary to Schubring's claim, Laugwitz did not attribute 20th century
number systems to Cauchy.  Rather, Laugwitz sought to understand
Cauchy's inferential moves in terms of their modern proxies.  May we
suggest that Schubring's mocking misrepresentation of Laugwitz's
position is, again, ``contrary to a style appropriate for a reasoned
scientific discussion.''  \cite[Section~5 `The Climax']{Sc15}

\subsection
{What charge is Laugwitz indicted on exactly?}
\label{s55}

In the abstract of his 1987 article in \emph{Historia Mathematica},
Laugwitz is careful to note that he interprets Cauchy's sum theorem
``with his [i.e., Cauchy's] own concepts":
\begin{quote}
It is shown that the famous so-called errors of Cauchy are correct
theorems when interpreted with his own concepts. \cite[p.~258]{La87}
\end{quote}
In the same abstract, Laugwitz goes on to emphasize: 
\begin{quote}
\emph{No assumptions} on uniformity or on nonstandard numbers are
needed. (emphasis added)
\end{quote}
Indeed, in section 7 on pages 264--266, Laugwitz gives a lucid
discussion of the sum theorem in terms of Cauchy's infinitesimals,
with not a whiff of modern number systems.  In particular this section
does not mention the article \cite{SL58}.  In a final section 15
entitled ``Attempts toward theories of infinitesimals," Laugwitz
presents a rather general discussion, with no specific reference to
the sum theorem, of how one might formalize Cauchyan infinitesimals in
modern set-theoretic terms.  A reference to \cite{SL58} appears in
this final section only.  Thus, Laugwitz carefully distinguishes
between his analysis of Cauchy's procedures, on the one hand, and the
ontological issues of possible implementations of infinitesimals in a
set-theoretic context, on the other.

Alas, all of Laugwitz's precautions went for naught.  In 2008, he
became a target of damaging innuendo in the updated version of
\emph{The Dictionary of Scientific Biography}.  Here C.~Fraser writes
as follows in his article on Cauchy:
\begin{quote}
Laugwitz's thesis is that certain of Cauchy's results that were
criticized by later mathematicians are in fact valid \emph{if one is
willing to accept certain assumptions} about Cauchy's understanding
and use of infinitesimals.  These assumptions reflect a theory of
analysis and infinitesimals that was worked out by Laugwitz and \ldots
Schmieden during the 1950s.%
\footnote{Fraser repeats the performance in 2015 when he claims that
``Laugwitz, \ldots{} some two decades following the publication by
Schmieden and him of the $\Omega$-calculus commenced to publish a
series of articles arguing that their non-Archimedean formulation of
analysis is well suited to interpret Cauchy's results on series and
integrals.''  \cite[p.~27]{Fr15} What Fraser fails to mention is that
Laugwitz specifically separated his analysis of Cauchy's
\emph{procedures} from attempts to account \emph{ontologically} for
Cauchy's infinitesimals in modern terms.}
\cite[p.~76]{Fr08} (emphasis added)
\end{quote}
Fraser and Schubring both claim that Laugwitz's interpretation of
Cauchy depends on assumptions that reflect a modern theory of
infinitesimals.  While Schubring charges Laugwitz with relying on
Robinson's hyperreals (see Section~\ref{s10}), Fraser's indictment is
based on the Omega-theory of Schmieden--Laugwitz.  Both Schubring and
Fraser are off the mark, as we showed above.

\section{Of infinitesimals and limits}

Schubring claims to distance himself from scholars like Grabiner, who
attribute to Cauchy a conceptual framework in the tradition of the
\emph{great triumvirate} of Cantor, Dedekind, and Weierstrass (see
\cite[p.~298]{Boy}):
\begin{quote}
I am criticizing historiographical approaches like that of Judith
Grabiner where one sees epsilon-delta already realized in Cauchy
\cite[Section~3]{Sc15}.
\end{quote}
Yet, when Schubring explains the so-called `compromise concept' that
he attributes to Cauchy, the infinitely small is governed by a
conception of limits in the context of ordinary (real) values.  While
Schubring is not guilty, like Grabiner, of ignoring the role of
infinitesimals in Cauchy's thinking, it is clear that Schubring
believes that it is the limit concept, rather than the infinitesimal
concept that is primary in this relationship.  A more detailed
analysis of such a Cauchy--Weierstrass tale appears in \cite{BK}.

Schubring's section (6.4), discussing the connection between limits
and the infinitely small, reveals that Schubring is a cheerleader for
the triumvirate team.  Thus, he writes that the infinitely small is
``\emph{subjugated to} \ldots the limit concept'' \cite[p.~454]{Sc05}
(emphasis added); furthermore, ``the \emph{infiniment petits}
represent only a subconcept of the \emph{limite} concept in Cauchy's
textbooks'' \cite[p.~455]{Sc05}.  Both claims seem to be based on
Cauchy's remark that ``a variable of this type has zero as a limit.''
However, this comment in Cauchy is not truly a part of the definition
of an infinitesimal, but rather a consequence of the definition which
precedes it.

\section{Indivisibles banned by the Jesuits}
\label{ban}

Most scholars agree that infinitesimal analysis was a natural
outgrowth of the indivisibilist techniques as developed by Galileo and
Cavalieri, and in fact Galileo's work may be closer to the
infinitesimal techniques of their contemporary Kepler.%
\footnote{Note that the term \emph{infinitesimal} itself was not
coined until the 1670s, by either Mercator or Leibniz; see
\cite[p.~63]{Le99}.}
Indivisibles were perceived as a theological threat and opposed on
doctrinal grounds in the 17th century \cite{Fe03}.  The opposition was
spearheaded by clerics and more specifically by the Jesuits.  Tracts
opposing indivisibles were composed by Jesuits Paul Guldin, Mario
Bettini, and Andr\'e Tacquet \cite[p.~291]{Re87}.  P.~Mancosu writes:
\begin{quote}
Guldin is taking Cavalieri to be composing the continuum out of
indivisibles, a position rejected by the Aristotelian orthodoxy as
having atomistic implications. \ldots{} Against Cavalieri's
proposition that ``all the lines" and ``all the planes" are magnitudes
- they admit of ratios - Guldin argues that ``all the lines \ldots{}
of both figures are infinite; but an infinite has no proportion or
ratio to another infinite."  \cite[p.~54]{Ma96}
\end{quote}
Tacquet for his part declared that the method of indivisibles ``makes
war upon geometry to such an extent, that if it is not to destroy it,
it must itself be destroyed.'' \cite{Al14}

In 1632 (the year Galileo was summoned to stand trial over
heliocentrism) the Society's Revisors General led by Jacob Bidermann
banned teaching indivisibles in their colleges \cite{Fe90},
\cite[p.~198]{Fe92}.  Referring to this ban, Feingold notes:
\begin{quote}
Six months later, General Vitelleschi formulated his strong opposition
to mathematical atomism in a letter he dispatched to Ignace Cappon in
Dole: ``As regards the opinion on \emph{quantity made up of
indivisibles}, I have already written to the Provinces many times that
it is in no way approved by me and up to now I have allowed nobody to
propose it or defend it. If it has ever been explained or defended, it
was done without my knowledge. Rather, I demonstrated clearly to
Cardinal Giovanni de Lugo himself that I did not wish our members to
treat or disseminate that opinion.''  \cite[p.~28-29]{Fe03} (emphasis
added)
\end{quote}
Indivisibles were placed on the Society's list of \emph{permanently}
banned doctrines in 1651 \cite{He96}.  In the 18th century, most
Jesuit mathematicians adhered to the methods of Euclidean geometry (to
the exclusion of the new infinitesimal methods):
\begin{quote}
\ldots le grand nombre des math\'ematiciens de [l'Ordre] resta
jusqu'\`a la fin du XVIII$^e$ si\`ecle profond\'ement attach\'e aux
m\'ethodes euclidiennes.  \cite[p.~77]{Bo27}
\end{quote}
Echoes of such bans were still heard in the 19th century.  Thus, in
1840, Moigno wrote:
\begin{quote}
In effect, either these magnitudes, smaller than any \emph{given}
magnitude, still have substance and are \emph{divisible}, or they are
simple and \emph{indivisible}: in the first case their existence is a
chimera,%
\footnote{Moigno's \emph{chimerical} anti-infinitesimal thread has not
remained without modern French adherents; see \cite{KKM}.}
since, necessarily greater than \emph{their half, their quarter},
etc., they are not actually less than any \emph{given} magnitude; in
the second hypothesis, they are no longer mathematical magnitudes, but
take on this quality, this would renounce the idea of the continuum
divisible to infinity, a necessary and fundamental point of departure
of all the mathematical sciences \ldots{} (as quoted in
\cite[p.~456]{Sc05}) [emphasis added]
\end{quote}
It is to be noted that Moigno formulates his objections to
infinitesimals specifically in the terminology of
indivisible/divisible which had been obsolete for nearly two
centuries.  Moigno saw a contradiction where there is none.  Indeed,
an infinitesimal is smaller than any \emph{assignable}, or
\emph{given} magnitude; and has been so at least since Leibniz.%
\footnote{Leibniz's dichotomy between assignable and inassignable
quantity, on which his concept of infinitesimal was based, finds a
rigorous mathematical treatment in the hyperreal number system (where
an assignable number is a standard real number).  Yet the mathematics
of earlier times allowed an adequate intuitive understanding of the
issue, sufficient to effectively and fruitfully use infinitesimals in
mathematical practice, even though a semantic base (accounting for the
ontology of a number) acceptable by modern standards was as yet
unavailable.  For further details on Leibniz's theoretical strategy in
dealing with infinitesimals see \cite{KS2}, \cite{KS1}, \cite{SK}.}
Thus, infinitesimals need not be less than ``their half, their
quarter, etc.,'' because the latter are not given/assignable.

In fact, Moigno's dichotomy is reminiscient of the dichotomy contained
in a critique of Galileo's indivisibles penned by Moigno's fellow
Jesuit Orazio Grassi some three centuries earlier.  P.~Redondi
summarizes it as follows:
\begin{quote}
As for light - composed according to Galileo of indivisible atoms,
more mathematical than physical - in this case, logical contradictions
arise.  Such indivisible atoms must be finite or infinite.  if they
are finite, mathematical difficulties would arise.  If they are
infinite, one runs into all the paradoxes of the separation to
infinity which had already caused Aristotle to discard the atomist
theory \ldots{} \cite[p.~196]{Re87}.
\end{quote}
This criticism appeared in the first edition of Grassi's book
\emph{Ratio ponderum librae et simbellae}, published in Paris in 1626.
According to Redondi, this criticism of Grassi's
\begin{quote}
exhumed a discounted argument, copied word-for-word from almost any
scholastic philosophy textbook. \ldots{} The Jesuit mathematician
[Paul] Guldin, great opponent of the geometry of indivisibles, and an
excellent Roman friend of [Orazio] Grassi, must have dissuaded him
from repeating such obvious objections.  Thus the second edition of
the \emph{Ratio}, the Neapolitan edition of 1627, omitted as
superfluous the whole section on indivisibles.  \cite[p.~197]{Re87}.
\end{quote}
Alas, unlike Father Grassi, Father Moigno had no Paul Guldin to
dissuade him.  Schubring fails to take into account the pertinent
Jesuit background of Moigno's rhetorical flourishes against
infinitesimals.  Schubring's endorsement (see Section~\ref{43}) of
Moigno's myopic anti-infinitesimal stance is nothing short of comical.

While we agree with Schubring on the importance of considering
context, we observe that he fails to consider the Jesuit context of
Moigno's flawed critique.  Schubring seeks to exploit Moigno's
\emph{text} in trying to prove a point about Cauchy.  However,
escaping the need of a genuine hermeneutic effort, Schubring ignores
the \emph{context} of Moigno's still being a Jesuit at the time he
wrote his 1840 \emph{text}.

\section{Gilain on records from the Ecole}

Cauchy wrote in the introduction to his 1821 \emph{Cours d'Analyse}:
\begin{quote}
In speaking of the continuity of functions, I could not dispense with
a treatment of the principal properties of infinitely small
quantities, properties which serve as the foundation of the
infinitesimal calculus. (translation from \cite{BS})
\end{quote}
When Cauchy writes that he was unable to dispense with infinitesimals,
was he complaining about being forced to teach infinitesimals, or was
he emphasizing the crucial importance of infinitesimals?
\cite[p.~436]{Sc05} asserts that Cauchy is complaining.  Schubring
bases his claim on the \emph{context} of Gilain's analysis of teaching
records from the \emph{Ecole} in \cite{Gi89}.  However, Gilain makes
the following points:
\begin{enumerate}
\item
Unlike Cauchy's later textbooks, his 1821 book was not commissioned by
the \emph{Ecole} but was rather written upon the personal request of
Laplace and Poisson.
\item
When the portion of the curriculum devoted to \emph{Analyse
Alg\'ebrique} was reduced in 1825, Cauchy insisted on placing the
topic of continuous functions [and therefore also the topic of
infinitesimals exploited in Cauchy's definition of continuity] at the
beginning of the Differential Calculus.
\end{enumerate}
Thus, Gilain writes:
\begin{quote}
\ldots rien ne montre, dans les documents utilis\'es, qu'il y ait eu
une hostilit\'e de principe de la part de Cauchy \`a la suppression de
l'analyse alg\'ebrique, en tant que partie autonome situ\'ee au
d\'ebut du cours d'analyse.  Ce qui lui importait, par contre,
c'\'etait la pr\'esence de plusieurs articles de cette partie, et les
m\'ethodes utilis\'ees pour les pr\'esenter. \ldots{} Notons, sur deux
points importants pour Cauchy, que les fonctions continues sont
plac\'ees au d\'ebut du calcul diff\'erentiel et que l'\'etude de la
convergence des s\'eries trouve sa place au voisinage de la formule de
Taylor, dans le calcul diff\'erentiel et int\'egral.%
\footnote{Translation: ``Nothing indicates, in the documents used, any
principled hostility on the part of Cauchy to the elimination of
algebraic analysis as an autonomous part placed at the beginning of
the analysis course.  What was important to him, on the other hand,
was the presence of several items from this part, and the methods used
in presenting them. \ldots{} Note two particular points important to
Cauchy, namely that continuous functions be placed at the beginning of
differential calculus, and that the study of the convergence of series
should find its place in the vicinity of Taylor's formula, in
differential and integral calculus.''}
\cite[end of \S\,131]{Gi89}
\end{quote}

The fact that Cauchy insisted on retaining continuity/infinitesimals
in 1825, as Gilain documents, indicates that the importance of
infinitesimals to Cauchy goes beyond their pedagogical value.  Also,
the fact that the \emph{Cours d'Analyse} was written on personal
request from Poisson and others, rather than being commissioned by the
\emph{Ecole}, indicates that the pressures of the type Schubring
claims are unlikely to have played a role in the writing of this
particular book.

Later texts were indeed commissioned by the \emph{Ecole}, but in 1821
Cauchy was free to write as he felt, and he felt that infinitesimals
were of crucial importance for analysis.  In painting a dim picture of
Cauchy's attitude toward infinitesimals, Schubring misjudges the
historical \emph{context} of the \emph{Cours d'Analyse} and ignores
the fact that decades after finishing his teaching stint at the
\emph{Ecole}, Cauchy exploits infinitesimals in an essential way in
his research article \cite{Ca53}.

\section{Conclusion}

Shoddy scholarship sometimes parades under the dual banner of
\emph{sophisticated historiographic analysis} and \emph{genuine
hermeneutical effort}.  As we have shown, Schubring's ``sophisticated
historiographical analysis" \cite[Section 5 `The Climax']{Sc15}
collapses in the face of simple historical facts and straightforward
mathematical arguments.  Moreover, Schubring's repeated
misrepresentation of Laugwitz's position is indeed ``contrary to a
style appropriate for a reasoned scientific discussion.''  (ibid.)

\section*{Acknowledgments}

The work of V. Kanovei was partially supported by RFBR grant
13-01-00006.  M.~Katz was partially funded by the Israel Science
Foundation grant no.~1517/12.  We are grateful to the anonymous
referees and to A. Alexander, R. Ely, and S. Kutateladze for their
helpful comments.  The influence of Hilton Kramer (1928-2012) is
obvious.

\medskip

\textbf{Piotr B\l aszczyk} is Professor at the Institute of
Mathematics, Pedagogical University (Cracow, Poland).  He obtained
degrees in mathematics (1986) and philosophy (1994) from Jagiellonian
University (Cracow, Poland), and a PhD in ontology (2002) from
Jagiellonian University.  He authored \emph{Philosophical Analysis of
Richard Dedekind's memoir} Stetigkeit und irrationale Zahlen (2008,
Habilitationsschrift).  He co-authored Euclid, Elements, Books V-VI.
Translation and commentary, 2013.  His research interest is in the
idea of continuum and continuity from Euclid to modern times.

\medskip

\textbf{Vladimir Kanovei} graduated in 1973 from Moscow State
University, and obtained a Ph.D. in physics and mathematics from
Moscow State University in 1976. In 1986, he became Doctor of Science
in physics and mathematics at Moscow Steklov Mathematical Institute
(MIAN).  He is currently Principal Researcher at the Institute for
Information Transmission Problems (IITP), Moscow, Russia. Among his
publications is the book \emph{Borel equivalence relations. Structure
and classification}. University Lecture Series 44. American
Mathematical Society, Providence, RI, 2008.

\medskip

\textbf{Mikhail G. Katz} (BA Harvard '80; PhD Columbia '84) is
Professor of Mathematics at Bar Ilan University, Ramat Gan, Israel.
He is interested in Riemannian geometry, infinitesimals, debunking
mathematical history written by the victors, as well as in true
infinitesimal differential geometry; see
\url{http://arxiv.org/abs/1405.0984}

\medskip

\textbf{David Sherry} is Professor of Philosophy at Northern Arizona
University, in the tall, cool pines of the Colorado Plateau.  He has
research interests in philosophy of mathematics, especially applied
mathematics and non-standard analysis.  Recent publications include
``Fields and the Intelligibility of Contact Action,'' \emph{Philosophy
90} (2015), 457--478.  ``Leibniz's Infinitesimals: Their Fictionality,
their Modern Implementations, and their Foes from Berkeley to Russell
and Beyond,'' with Mikhail Katz, \emph{Erkenntnis 78} (2013), 571-625.
``Infinitesimals, Imaginaries, Ideals, and Fictions,'' with Mikhail
Katz, \emph{Studia Leibnitiana 44} (2012), 166--192.  ``Thermoscopes,
Thermometers, and the Foundations of Measurement,'' \emph{Studies in
History and Philosophy of Science 24} (2011), 509--524.  ``Reason,
Habit, and Applied Mathematics,'' \emph{Hume Studies 35} (2009),
57-85.

\end{document}